\documentclass[12pt,centertags]{amsart}
\usepackage{amsmath,amstext,amsthm,a4,amssymb}
\numberwithin{equation}{section}

\newcommand{\field}[1]{\mathbb{#1}}
\newcommand{\Z}{\field{Z}}
\newcommand{\R}{\field{R}}
\newcommand{\C}{\field{C}}
\newcommand{\N}{\field{N}}

\newcommand{\til}[1]{\widetilde{#1}}
\newcommand{\tx}{\til{X}}
\newcommand{\tj}{\til{J}}

\newcommand{\tsx}{\til{x}}
\newcommand{\ty}{\til{y}}
\newcommand{\tom}{\til{\omega}}
\newcommand{\tl}{\til{L}}
\newcommand{\te}{\til{E}}
\newcommand{\tdk}{\til{D}_k}
\newcommand{\tdel}{\til{\Delta}_k}
\newcommand{\tp}{\til{P}}
\newcommand{\tf}{\til{F}}

\newcommand{\tsc}{\til{\Delta}^\#_k}
\newcommand{\cali}[1]{\mathcal{#1}}
\newcommand{\Ha}{\cali{H}}

\DeclareMathOperator{\End}{End}
\DeclareMathOperator{\ke}{Ker}
\DeclareMathOperator{\Dom}{Dom}
\DeclareMathOperator{\ran}{Range}
\DeclareMathOperator{\rank}{rank}
\DeclareMathOperator{\Id}{Id}

\DeclareMathOperator{\tr}{Tr}
\DeclareMathOperator{\ind}{index}
\DeclareMathOperator{\td}{Td}
\DeclareMathOperator{\vol}{vol}
\DeclareMathOperator{\ch}{ch}
\newcommand{\spin}{$\text{spin}^c$ }

\newcommand{\norm}[1]{\lVert#1\rVert}

\newcommand{\om}{\omega} 
\newcommand{\g}{\Gamma}

\newtheorem{theorem}{Theorem}[section]
\newtheorem{lemma}[theorem]{Lemma}
\newtheorem{proposition}[theorem]{Proposition}
\newtheorem{corollary}[theorem]{Corollary}

\theoremstyle{definition} 
\newtheorem{defn}[theorem]{Definition}
\theoremstyle{remark}
\newtheorem{rem}[theorem]{Remark}

\newcommand{\comment}[1]{}

\begin{document}
\title{The $\text{SPIN}^c$ Dirac operator on high tensor powers of a line bundle}
\author{Xiaonan Ma  and George Marinescu}
\thanks{Humboldt-Universit\"at zu Berlin, Institut f\"ur Mathematik, Rudower Chaussee 25,
12489 Berlin, Germany. (xiaonan@mathematik.hu-berlin.de, george@mathematik.hu-berlin.de). 
Partially supported by SFB 288.} 
\address{}
\date{}
\begin{abstract}
We study the asymptotic of the spectrum of the \spin Dirac operator on high tensor powers of a line bundle.
As application, we get a simple proof of the main result of Guillemin--Uribe \cite[Theorem 2]{GU},
which was originally proved by using the analysis of Toeplitz operators of Boutet de Monvel and Guillemin
\cite{BG}.
\end{abstract} 

\maketitle

\section{Introduction} \label{s1}

Let $(X,\om)$ be a compact symplectic manifold of real dimension $2n$.
Assume that there exists a hermitian line bundle $L$ over $X$ endowed with
a hermitian connection $\nabla^L$ with the property that
$\frac{\sqrt{-1}}{2\pi}R^L=\omega$, where $R^L=(\nabla^L)^2$ is the curvature of $(L,\nabla^L)$.
Let $E$ be a hermitian vector
bundle $E$ on $X$.

Let $g^{TX}$ be a riemannian metric on $X$.
Let $J_0:TX\longrightarrow TX$ be the skew--adjoint linear map which satisfies
the relation $\om(u,v)=g^{TX}(J_{0}u,v)$. Let $J$ be an almost complex structure
which is compatible with $g^{TX}$ and $\om$. 
Then one can construct canonically a
\spin  Dirac operator $D_k$ acting on 
$\Omega^{0,{\scriptscriptstyle{\bullet}}}(X,L^k\otimes E)=\bigoplus_{q=0}^n\Omega^{0,q}(X,L^k\otimes E)$,
the direct sum of spaces of $(0,q)$--forms with values in $L^k\otimes E$.
Let $\lambda=\displaystyle\inf_{{u\in T_x^{(1,0)}X,\,x\in X}} R^L_x(u,\overline{u})/|u|^2_{g^{TX}} >0$\,.
One of our main results is as follows: 
\begin{theorem}\label{specDirac}
There exists $C>0$ such that for $k\in \N$,  the spectrum of $D^2_k$ is contained in the
set $\{0\}\cup(2k\lambda-C,+\infty)$. Set $D_k^{-}=D_k\upharpoonright_{\Omega^{0,{\text{odd}}}}$, 
then for $k$ large enough,
we have 
\begin{equation}\label{vanish}
\ke D_k^{-}=\{0\}.
\end{equation}
\end{theorem}
\noindent
We recover with \eqref{vanish} the vanishing result of
\cite[Theorem 2.3]{BU}, \cite[Theorem 3.2]{Bra1}.
Another interesting application is to describe the asymptotics of the 
spectrum of the metric Laplacian 
$\Delta_k=(\nabla^{L^k\otimes E})^{\ast}\,\nabla^{L^k\otimes E}$ acting on
$\mathcal{C}^{\infty}(X,L^k\otimes E)$.
Introduce the smooth function
$\tau(x)=\sum_j R^L (w_j,\overline{w}_j) >0$, $x\in X$, where
$\{w_j\}_{j=1}^n$
is an orthonormal basis of $T_x^{(1,0)}X$. 
\begin{corollary}\label{schrodinger}
The spectrum of the Schr\"odinger
operator $\Delta^\#_k=\Delta_k-k\tau$
is contained in the union $(-a,a)\cup(2k\lambda-b,+\infty)$, 
where $a$ and $b$ are positive constants independent of $k$. 
For $k$ large enough, the number $d_k$ of eigenvalues on the interval $(-a,a)$ satisfies
$d_k=\langle\ch(L^k\otimes E)\td(X),[X]\rangle$. 
In particular $d_k\sim k^n(\rank{E})\vol_{\om}(X)$.
\end{corollary}
In the case $E$ is a trivial line bundle, Corollary \ref{schrodinger} is the main result of Guillemin and Uribe 
\cite[Theorem 2]{GU} \footnote{In \cite{GU}, they only knew  $d_k \sim k^n \vol_{\om}(X)$.  When $J_0=J$,  
Borthwick and Uribe \cite[p. 854]{BU}  got 
the precise value $d_k$, for large enough $k$, in this case.}. The idea in \cite{BU}, \cite{Bra1}, \cite{Bra}, \cite{GU} is that one first reduces the problem
to a problem on the unitary circle bundle of $L^\ast$, then one applies Melin inequality \cite[Theorem 22.3.2]
{Hor} to show that $\Delta^\#_k$ is semi--bounded from below. In order to prove \cite[Theorem 2]{GU},
they  apply the analysis of Toeplitz structures
of Boutet de Monvel--Guillemin \cite{BG}.
For the interesting applications of \cite[Theorem 2]{GU}, we refer the reader to Borthwick and Uribe
\cite{BU}, \cite{BU1}, \cite{BU2}. 
For the related topic on geometry quantization, see \cite{M}, \cite{TZ}.
Our proof is based on a direct application of Lichnerowicz formula.

This paper is organized as follows. In Section \ref{s2}, we recall the construction of the \spin Dirac operator
and prove our  main technical result, Theorem \ref{main}. In  Section \ref{s3}, we  prove  Theorem \ref{specDirac}
and  Corollary \ref{schrodinger}. In Section \ref{s4}, we generalize our result to the $L_2$ case.
In particular, we obtain a new proof of  \cite[Theorem 2.6]{Bra}. 

\section{The Lichnerowicz formula}\label{s2}

Let $(X,\om)$ be a compact symplectic manifold. Let $(L,h^L)$, $(E,h^E)$ be two hermitian
vector bundles  endowed with hermitian connections 
$\nabla^L$ and $\nabla^E$ respectively. Let $R^L$ and $R^E$ be their curvatures. 
We assume $\rank L=1$ and $R^L=-2\pi\sqrt{-1}\om$.
Let $g^{TX}$ be an arbitrary riemannian metric on $TX$.  
Let $J$ be an almost complex structure  which is
compatible with $g^{TX}$ and $\om$ 
(For the existence of $J$, we refer to \cite[p.61]{MS}).  
Then $J$ defines canonically an orientation of $X$. 
Let  $J_0:TX\longrightarrow TX$ be  the skew--adjoint linear map defined by
\begin{equation}\label{b1}
\om(u,v)=g^{TX}(J_{0}u,v), \ \ {\rm for} \ \ u,v \in TX.
\end{equation}
Then  $J$  commutes with $J_0$.

Let $TX^c=TX\otimes_\R \C$ denote the complexification of the 
tangent bundle. The almost complex structure $J$ induces a splitting
$TX^c=T^{(1,0)}X\oplus T^{(0,1)}X$, where $T^{(1,0)}X$ and $T^{(0,1)}X$
are the eigenbundles of $J$ corresponding to the eigenvalues $\sqrt{-1}$
and $-\sqrt{-1}$ respectively. Accordingly, we have a decomposition of the
complexified cotangent bundle: $T^{\ast}X^c=T^{(1,0)\,\ast}X\oplus
T^{(0,1)\,\ast}X$. The exterior algebra bundle decomposes as
$\Lambda T^{\ast}X^c=\oplus_{p,q}\Lambda^{p,q}$, where
$\Lambda^{p,q}:=\Lambda^{p,q}T^{\ast}X^c=\Lambda^{p}(T^{(1,0)\,\ast}X)\otimes
\Lambda^{q}(T^{(0,1)\,\ast}X)$.

Let $\nabla^{TX}$ be the Levi--Civita connection
of the metric $g^{TX}$, and let $\nabla^{1,0}$ and $\nabla^{0,1}$ be the canonical 
hermitian connections on $T^{(1,0)}X$ and $T^{(0,1)}X$ respectively: 
\begin{align*}
\nabla^{1,0}&=\tfrac{1}{4}(1-\sqrt{-1} J)\,\nabla^{TX}\,(1-\sqrt{-1} J)\,,\\
\nabla^{0,1}&=\tfrac{1}{4}(1+\sqrt{-1} J)\,\nabla^{TX}\,(1+\sqrt{-1} J)\,.
\end{align*}
Set $A_2=\nabla^{TX}-\big(\nabla^{1,0}\oplus\nabla^{0,1}\big)$  $\in T^{\ast}X\otimes\End(TX)$ which satisfies $J\,A_2=-A_2\,J$.

Let us recall some basic facts about the \spin Dirac operator on an almost complex 
manifold \cite[Appendix D]{LM}. 
The fundamental $\Z_2$
spinor bundle induced by $J$ is given by 
$\Lambda^{0,\scriptscriptstyle\bullet}=\Lambda^{\text{even}}
(T^{(0,1)\,\ast}X)\oplus\Lambda^{\text{odd}} (T^{(0,1)\,\ast}X)$.
For any $v\in TX$ with decomposition $v=v_{1,0}+v_{0,1} 
\in T^{(1,0)}X\oplus T^{(0,1)}X$,  let ${\overline v^\ast_{1,0}}\in T^{(0,1)*}X$ be the metric 
dual of $v_{1,0}$. Then $\mathbf{c}(v)=\sqrt{2}({\overline v^\ast_{1,0}}\wedge-
i_{v_{\,0,1}})$ defines the Clifford action of $v$ on 
$\Lambda^{0,\scriptscriptstyle\bullet}$, where $\wedge$ and $i$ 
denote the exterior and interior product respectively.

Formally, we may think
$$ \Lambda^{0,\scriptscriptstyle\bullet}=S\left(TX\right)\otimes 
\left(\det{T^{(1,0)}X}\right)^{1/2}\,,$$ 
where $S\left(TX\right)$ is the spinor bundle of the possibly non--existent spin 
structure on $TX$, and $\left(\det{T^{(1,0)}X}\right)^{1/2}$ is the possibly 
non--existent square root of $\det{T^{(1,0)}X}$. 

Moreover, by \cite[pp.397--398]{LM}, $\nabla^{TX}$ induces canonically a Clifford connection 
on $\Lambda^{0,\scriptscriptstyle{\bullet}}$.
Formally, let $\nabla^{S(TX)}$ be the Clifford connection on $S(TX)$ induced by $\nabla^{TX}$,
and let $\nabla^{\det}$ be the connection on $(\det{T^{(1,0)}X})^{1/2}$ induced by
$\nabla^{1,0}$. Then
$$
\nabla^{\text{Cliff}}=\nabla^{S(TX)}\otimes\Id+\Id\otimes\nabla^{\det}\,.
$$ 
 Let $\{w_j\}_{j=1}^n$
be a local orthonormal frame of $T^{(1,0)}X$. Then
\begin{equation}\label{frame1}
e_{2j}=\tfrac{1}{\sqrt{2}}(w_j+\overline{w}_j)\quad\text{and}\quad
e_{2j-1}=\tfrac{\sqrt{-1}}{\sqrt{2}}(w_j-\overline{w}_j)\,, \quad j=1,\dotsc,n\,,
\end{equation}
form an orthonormal frame of $TX$.
Let $\{w^j\}_{j=1}^n$ be the dual frame of $\{w_j\}_{j=1}^n$.
Let $\Gamma$ be the connection form of $\nabla^{1,0}\oplus \nabla^{0,1}$ in local coordinates.
Then $\nabla^{TX}=d+\Gamma+A_2$. By \cite[Theorem 4.14, p.110]{LM},  the Clifford connection
$\nabla^{\text{Cliff}}$ on $\Lambda^{0,\scriptscriptstyle{\bullet}}$
has the following local form:
\begin{equation}\label{local}
\begin{split}
\nabla^{\text{Cliff}}&= d + \tfrac{1}{4} \sum_{i,j} \big\langle (\Gamma+ A_2)e_i, e_j\big\rangle \mathbf{c}(e_i)\mathbf{c}(e_j)\\
&=d+\sum_{l,m}\Big\lbrace\big\langle\Gamma w_l,\overline{w}_m
\big\rangle\,\overline{w}^l\wedge\,i_{\overline{w}_m}+\\
&\hspace*{7mm}\tfrac{1}{2}\big\langle A_2 w_l,w_m
\big\rangle\,i_{\overline{w}_l}\,i_{\overline{w}_m}+
\tfrac{1}{2}\big\langle A_2\overline{w}_l,\overline{w}_m
\big\rangle\,\overline{w}^l\wedge\,\overline{w}^m\wedge \Big\rbrace\,.
\end{split}
\end{equation}

Let $\nabla^{L^k\otimes E}$ be  the connection on $L^k \otimes E$ induced 
by $\nabla^L, \nabla^E$. Let $\nabla^{\Lambda^{0,\scriptscriptstyle{\bullet}}\otimes L^k\otimes E}$ be the connection on $\Lambda^{0,\scriptscriptstyle{\bullet}}\otimes L^k\otimes E$,
\begin{equation}
\nabla^{\Lambda^{0,\scriptscriptstyle{\bullet}}\otimes L^k\otimes E}=
\nabla^{\text{Cliff}}\otimes\Id+\Id\otimes\nabla^{L^k\otimes E}.
\end{equation}

Along the fibers of $\Lambda^{0,\scriptscriptstyle{\bullet}}\otimes L^k\otimes E$,
we consider the pointwise scalar product $\langle\cdot,\cdot\rangle$  
induced by $g^{TX}$, $h^L$ and $h^E$. 
Let $dv_X$ be the riemannian volume form of $(TX, g^{TX})$.
The $L_2$--scalar product on $\Omega^{0,\scriptscriptstyle{\bullet}}(X,L^k\otimes E)$, the space of smooth sections of 
$\Lambda^{0,\scriptscriptstyle{\bullet}}\otimes L^k\otimes E$,
is given by 
\begin{equation}\label{l2}
(s_1,s_2)=\int_X\langle s_1(x),s_2(x)\rangle\,dv_X(x)\,.
\end{equation}
We denote the corresponding norm with $\norm{\cdot}$. 

\begin{defn}\label{Dirac}
The \spin Dirac operator $D_k$ is defined by 
\begin{equation}\label{defDirac}
D_k=\sum_{j=1}^{2n}\mathbf{c}(e_j)\nabla^{\Lambda^{0,\scriptscriptstyle{\bullet}}\otimes L^k\otimes E}_{e_j}:
\Omega^{0,\scriptscriptstyle{\bullet}}(X,L^k\otimes E)\longrightarrow
\Omega^{0,\scriptscriptstyle{\bullet}}(X,L^k\otimes E)\,.
\end{equation}
\end{defn}
\noindent
$D_k$ is a formally self--adjoint, first order elliptic differential operator on 
$\Omega^{0,\scriptscriptstyle{\bullet}}(X,L^k\otimes{E})$,
which interchanges $\Omega^{0,\text{even}}(X,L^k\otimes E)$
and $\Omega^{0,\text{odd}}(X,L^k\otimes E)$. We denote
\begin{equation}\label{}
D_k^{+}=D_k\upharpoonright_{\Omega^{0,{\text{even}}}}, \quad D_k^{-}=D_k\upharpoonright_{\Omega^{0,{\text{odd}}}}.
\end{equation}

Let $R^{T^{(1,0)}X}$ be the curvature of $\big(T^{(1,0)}X,\nabla^{1,0}\big)$.
Let 
\begin{align}
&\om_d=-\sum_{l,m} R^L (w_l,\overline{w}_m)\,\overline{w}^m\wedge
\,i_{\overline{w}_l}\,,\\
&\tau(x)=\sum_j R^L (w_j,\overline{w}_j)\,.\nonumber
\end{align}
Remark that by (\ref{b1}), at $x\in X$, there exists $\{ w_i\}_{i=1}^n$ an orthogonal basis of $T^{(1,0)} X$, 
such that $J_0= \sqrt{-1}\  {\rm diag} (a_1 (x), \cdots, a_n(x))\in {\rm End} (T^{(1,0)} X)$, and $a_i(x)>0$ 
for $i\in \{ 1, \cdots, n\}$. So
\begin{align}
&\om_d=- 2 \pi \sum_{l} a_l(x)\,\overline{w}^l\wedge
\,i_{\overline{w}_l}\,,\\
&\tau(x)=2 \pi \sum_l a_l(x)\,.\nonumber
\end{align}
The following Lichnerowicz formula is crucial for us.
\begin{theorem} \label{Lichn}
The square of the Dirac operator satisfies the equation:
\begin{equation}\label{Lich}
D^2_k=\left(\nabla^{\Lambda^{0,\scriptscriptstyle{\bullet}}\otimes L^k\otimes E}\right)^{\ast}\,
\nabla^{\Lambda^{0,\scriptscriptstyle{\bullet}}\otimes L^k\otimes E}-2k\om_d-
k\tau+\tfrac{1}{4}K+ \mathbf{c}(R),
\end{equation}
where $K$ is the scalar curvature of $(TX,g^{TX})$, and
\begin{equation*}
\mathbf{c}(R)=\sum_{l<m}\left(R^E+\tfrac{1}{2}\tr\left[R^{T^{(1,0)}X}\right]
\right)(e_l,e_m)\,\mathbf{c}(e_l)\,\mathbf{c}(e_m)\,.
\end{equation*}
\end{theorem} 
\begin{proof}
By Lichnerowicz formula \cite[Theorem 3.52]{BVG}, we know that
\begin{equation}\label{}
D^2_k=\left(\nabla^{\Lambda^{0,\scriptscriptstyle{\bullet}}\otimes L^k\otimes E}\right)^{\ast}\,\nabla^{\Lambda^{0,\scriptscriptstyle{\bullet}}\otimes L^k\otimes E}+\tfrac{1}{4}K+ 
\mathbf{c}(R)+k\sum_{l<m}R^L(e_l,e_m)\,\mathbf{c}(e_l)\,\mathbf{c}(e_m)\,.
\end{equation}
Now, we identify $R^L$ with a purely imaginary antisymmetric matrix $- 2 \pi \sqrt{-1} J_0\in {\rm End} (T X)$ by 
(\ref{b1}). As $J_0\in {\rm End} (T^{(1,0)} X)$, by   \cite[Lemma 3.29]{BVG}, we get (\ref{Lich}).
\end{proof}
\begin{rem}\label{rem2.3}
Let $\mathcal{E}=\mathcal{E}^{+}\oplus\mathcal{E}^{-}$ 
be a Clifford module. 
Then it was observed by Braverman \cite[\S 9]{Bra1} that, with the same proof of \cite[Proposition 3.35]{BVG}, 
there exists a vector bundle $W$ on $X$ such that 
$\mathcal{E}=\Lambda^{0,\scriptscriptstyle{\bullet}}\otimes W$ as a $\Z_2$--graded
Clifford module. This means (\ref{vanish}) in fact recovers
 \cite[Theorem 3.2]{Bra1}.
\end{rem}

As a simple consequence of Theorem \ref{Lichn}, we recover the statement
on the drift of spectrum of the metric Laplacian first proved by Guillemin--Uribe 
\cite[Theorem 1]{GU}, (see also \cite[Theorem 2.1]{BU}, \cite[Theorem 4.4]{Bra1}),
by passing to the circle bundle of $L^\ast$ and applying Melin's inequality
\cite[Theorems 22.3.2--3]{Hor}.
\begin{corollary}\label{bounded}
There exists $C>0$ such that for $k\in\N$, the metric Laplacian $\Delta_k=
\big(\nabla^{L^k\otimes E}\big)^{\ast}\,\nabla^{L^k\otimes E}$
on $\mathcal{C}^\infty(X,L^k\otimes E)$  satisfies \,:
\begin{equation}\label{est1}
\Delta_k-k\tau\geqslant-C\,.
\end{equation}
\end{corollary}
\begin{proof}
By (\ref{Lich}),  $s\in\mathcal{C}^\infty(X,L^k\otimes E)$, 
\begin{equation}\label{est2}
\norm{D_{k}s}^2=\norm{\nabla^{\Lambda^{0,\scriptscriptstyle{\bullet}}\otimes L^k\otimes E}s}^2-k(\tau(x)s,s)+ 
\left(\left(\tfrac{1}{4}K+\mathbf{c}(R)\right)s,s\right)\,.
\end{equation}
>From \eqref{local}, we infer that 
$$\big\|\nabla^{\Lambda^{0,\scriptscriptstyle{\bullet}}\otimes L^k\otimes E}s\big\|^2=
\big\|\nabla^{L^k\otimes E}s\big\|^2+
\Big\|\sum_{l,m}\big\langle A_2\overline{w}_l,\overline{w}_m
\big\rangle\,\overline{w}^l\wedge\overline{w}^m\wedge{s}\Big\|^2\,.$$
and therefore there exists a constant $C>0$ not depending on $k$ such that
$$ 
0\leqslant\norm{D_{k}s}^2\leqslant\big\|\nabla^{L^k\otimes E}s\big\|^2
-k(\tau(x)s,s)+C\norm{s}^2
=\big((\Delta_{k}-k\tau(x))s,s\big)+C\norm{s}^2\,.
$$
\end{proof}
The following is our main technical result.
\begin{theorem}\label{main}
There exists $C>0$ such that for any $k\in \N$ and any $s\in\Omega^{>0}(X,L^k\otimes E)
=\bigoplus_{q\geqslant 1}\Omega^{0,q}(X,L^k\otimes E)$, 
\begin{equation}\label{main1}
\norm{D_{k}s}^2\geqslant(2k\lambda-C)\norm{s}^2\,.
\end{equation}
\end{theorem}
\begin{proof}
By \eqref{Lich}, for $s\in\Omega^{0,\scriptscriptstyle{\bullet}}(X,L^k\otimes E)$\,,
\begin{equation}\label{lich1}
\norm{D_{k}s}^2=\lbrace\norm{\nabla^{\Lambda^{0,\scriptscriptstyle{\bullet}}\otimes L^k\otimes E}s}^2-k(\tau(x)s,s)\rbrace
-2k(\om_{d}s,s)+
\left(\left(\tfrac{1}{4}K+\mathbf{c}(R)\right)s,s\right)\,.
\end{equation}
We consider now $s \in \mathcal{C}^\infty(X,L^k\otimes E')$, where $E^\prime=E\otimes\Lambda^{0,\scriptscriptstyle{\bullet}}$.  
Estimate \eqref{est1} becomes
\begin{equation}\label{est1prim}
\big\|\nabla^{L^k\otimes E^\prime}s\big\|^2-k(\tau(x)s,s)\geqslant-C\norm{s}^2\,.
\end{equation}
If $s\in\Omega^{>0}(X,L^k\otimes E)$, the second term of (\ref{lich1}), $-2k(\om_{d}s,s)$ is bounded below by
$2k\lambda\norm{s}^2$. While the third term of (\ref{lich1}) is $O(\norm{s}^2)$.
The proof of  \eqref{main1} is completed.
\end{proof}
\section{Applications of Theorem \ref{main}}\label{s3}

\begin{proof}[Proof of Theorem \ref{specDirac}] By (\ref{main1}), we get immediately (\ref{vanish}). For the rest, 
we use the trick of the proof of Mckean--Singer formula.

Let  $\Ha_{\mu}$  be the spectral space of $D^2_k$ corresponding to  
the interval $(0,\mu)$. 
Let $\Ha^{+}_{\mu}$, $\Ha^{-}_{\mu}$ be 
the intersections of $\Ha_{\mu}$ with the spaces of forms of even and odd 
degree respectively. 
Then
$\Ha_{\mu}=\Ha^{+}_{\mu}\oplus\Ha^{-}_{\mu}$.
Since $D^{+}_k$ commutes with the spectral projection,
we have a well defined operator  
$D^{+}_k:\Ha^{+}_{\mu}\longrightarrow\Ha^{-}_{\mu}$ which is obviously  
injective. But estimate \eqref{main1} implies that $\Ha^{-}_{\mu}=0$
for every $\mu<2k\lambda-C$, hence also $\Ha^{+}_{\mu}=0$, for this range 
of $\mu$. Thus $\Ha_{\mu}=0$, for $0<\mu<2k\lambda-C$.
The proof of our theorem  is completed.
\end{proof}
\begin{proof}[Proof of Corollary \ref{schrodinger}]Let 
$P_k:\Omega^{0,\scriptscriptstyle{\bullet}}(X,L^k\otimes E)\longrightarrow
\mathcal{C}^{\infty}(X,L^k\otimes E)$ be the orthogonal projection.
For $s\in\Omega^{0,\scriptscriptstyle{\bullet}}(X,L^k\otimes E)$, we will denote
$s_0=P_{k}s$ its $0$\,--\,degree component.
We will estimate 
$\Delta^\#_k$ on $P_k(\ke D_k^{+})$ and $(\ke D_k^{+})^\perp\cap\mathcal{C}^{\infty}(X,L^k\otimes E)$.

In the sequel we denote with $C$ all positive constants independent of $k$, although there may be different constants
for different estimates. 
>From \eqref{est2}, there exists $C>0$ such that for $s\in\mathcal{C}^{\infty}(X,L^k\otimes E)$,
\begin{equation}\label{est4}
\big|\norm{D_{k}s}^2-(\Delta^\#_{k}s,s)\big|\leqslant C\norm{s}^2\,.
\end{equation}
Theorem \ref{specDirac} and  \eqref{est4} show that there exists  $b>0$ such that
for $k\in \N$, 
\begin{equation}\label{est5}
(\Delta^\#_{k}s,s)\geqslant(2k\lambda-b)\norm{s}^2\,,\quad {\rm for } \ s\in\mathcal{C}^{\infty}(X,L^k\otimes E)\cap(\ke D_k^{+})^\perp. 
\end{equation} 

We focus now on elements from $P_k(\ke D_k^{+})$,
and assume $s\in\ke D_k^{+}$. Set
$s^\prime=s-s_0\in\Omega^{>0}(X,L^k\otimes E)$.
By \eqref{lich1}, \eqref{est1prim},
\begin{equation}\label{bu0}
-2k(\om_d{s},s)\leqslant{C}\norm{s}^2\,.
\end{equation}
We obtain thus \cite[Theorem 2.3]{BU} (see also \cite{BU3}, \cite[Theorem 3.13]{Bra1}) 
for $k\gg1$, 
\begin{equation}\label{bu}
\norm{s^\prime}\leqslant{C}k^{-1/2}\norm{s_0}\,.
\end{equation}
(from (\ref{bu}), they got 
${\rm Ker} D^-_k =0$ for $k\gg1$). In view of \eqref{lich1} and \eqref{bu},
\begin{equation}\label{est8}
\norm{\nabla^{\Lambda^{0,\scriptscriptstyle{\bullet}}\otimes L^k\otimes E}s}^2
-k(\tau(x)s_0,s_0)\leqslant{C}\norm{s_0}^2\,.
\end{equation}
By \eqref{local}, 
\begin{equation}\label{local1}
\nabla^{\Lambda^{0,\scriptscriptstyle{\bullet}}\otimes L^k\otimes E}s=
\nabla^{L^k\otimes E}s_0+A'_{2}s_2+\alpha\,,
\end{equation}
where $s_2$ is the component of degree $2$ of $s$, $A'_2$ is a contraction operator comming 
from the middle term of \eqref{local}, and $\alpha\in\Omega^{>0}(X,L^k\otimes E)$.
By \eqref{est8}, \eqref{local1}, we know
\begin{equation}\label{est9}
\big\|\nabla^{L^k\otimes E}s_0+A'_{2}s_2\big\|^2-k(\tau(x)s_0,s_0)\leqslant{C}\norm{s_0}^2\,,
\end{equation}
and by \eqref{bu}, \eqref{est9},
\begin{equation}\label{est14}
\big\|\nabla^{L^k\otimes E}s_0\big\|^2\leqslant{C}k\norm{s_0}^2\,,
\end{equation}
By \eqref{bu} and \eqref{est14}, we get
\begin{equation}\label{est10}
\begin{split}
\big\|\nabla^{L^k\otimes E}s_0+A'_{2}s_2\big\|^2&\geqslant\big\|\nabla^{L^k\otimes E}s_0\big\|^2-
2\big\|\nabla^{L^k\otimes E}s_0\big\|\big\|A'_{2}s_2\big\|\\
&\geqslant\big\|\nabla^{L^k\otimes E}s_0\big\|^2-C\norm{s_0}^2.
\end{split}
\end{equation}
Thus, \eqref{est9} and \eqref{est10} yield
\begin{equation}\label{est12}
\big\|\nabla^{L^k\otimes E}s_0\big\|^2-k(\tau(x)s_0,s_0)\leqslant{C}\norm{s_0}^2\,.
\end{equation}
By \eqref{est1} and \eqref{est12}, there exists a constant $a>0$ such that 
\begin{equation}\label{est11}
\big|\big(\Delta^\#_{k}s,s\big)\big|\leqslant a\norm{s}^2\,,\quad s\in P_k(\ke D_k^{+})\,.
\end{equation}
By \eqref{bu}, we know that for $k\gg1$, $P_k:\ke D^{+}_k\longrightarrow P_k(\ke D^{+}_k)$
is bijective, and 
\begin{equation}\label{est15} 
\mathcal{C}^{\infty}(X,L^k\otimes E)=
P_k(\ke D_k^{+})\oplus(\ke D_k^{+})^\perp\cap\mathcal{C}^{\infty}(X,L^k\otimes E)\,.
\end{equation}

The proof is now reduced to a direct application of the minimax principle for
the  operator $\Delta^\#_{k}$. It is clear that \eqref{est5}
and \eqref{est11} still hold for elements in the Sobolev space $W^1(X,L^k\otimes E)$,
which is the domain of the quadratic form
$Q_k(f)=\big\|\nabla^{L^k\otimes E}f\big\|^2-k(\tau(x)f,f)$ associated to $\Delta^\#_{k}$.
Let  $\mu^k_1\leqslant\mu^k_2\leqslant\cdots\leqslant \mu^k_j\leqslant\cdots$ 
($j\in\N$) be the eigenvalues of $\Delta^\#_{k}$. 
Then, by the minimax principle \cite[pp.76--78]{RS},
\begin{equation}\label{min}
\mu^k_j=\min_{F\subset\Dom{Q_k}}\max_{f\in F\,,\,\norm{f}=1}Q_k(f)\,.
\end{equation} 
where $F$ runs over the subspaces of dimension $j$ of $\Dom{Q_k}$.

By \eqref{est11}
and \eqref{min}, we know  $\mu^k_j\leqslant a$, 
for $j\leqslant\dim\ke D_k^{+}$.
Moreover, any subspace $F\subset\Dom{Q_k}$ with $\dim{F}\geqslant\dim\ke D_k^{+}+1$ contains
an element $0\neq{f}\in{F}\cap(\ke D_k^{+})^\perp$. 
By \eqref{est5}, \eqref{min},
we obtain $\mu^k_j\geqslant 2 k\lambda-b$, for $j\geqslant\dim\ke D_k^{+}+1$.

By Theorem \ref{specDirac} and Atiyah--Singer theorem \cite{AS},
\begin{equation}\label{index}
\dim\ke D_k^{+}=\ind D_k^{+}=\langle\ch(L^k\otimes E)\td(X),[X]\rangle 
\end{equation}
where $\td(X)$ is the Todd class of an almost complex structure 
compatible with $\om$. The index is a polynomial in $k$ of degree 
$n$ and of leading term $k^n(\rank{E})\vol_{\om}(X)$, where $\vol_{\om}(X)$ 
is the symplectic volume of $X$.

The proof of our corollary is completed.
\end{proof}
\begin{rem}
If $(X,\om)$ is K\"ahler and if  $L$, $E$ are holomorphic vector bundles, then
$D_k=\sqrt{2}\big(\overline{\partial}+\overline{\partial}^{\,\ast}\big)$ where
$\overline{\partial}=\overline{\partial}^{\,L^k\otimes{E}}$.  $D^2_k$ preserves the
$\Z$--grading of $\Omega^{0,\scriptscriptstyle{\bullet}}$. By using the Bochner--Kodaira--Nakano
formula, Bismut and Vasserot \cite[Theorem 1.1]{BiV} proved Theorem \ref{main}. As
$\overline{\partial}:(\ke D_k^{+})^\perp\cap\mathcal{C}^{\infty}(X,L^k\otimes E)
\longrightarrow\Omega^{0,1}(X,L^k\otimes E)$ is injective, we infer
\begin{equation}
 2 \big\|\overline{\partial}s\big\|^2\geqslant(2k\lambda-C)\norm{s}^2,
\ {\rm for} \  s\in(\ke D_k^{+})^\perp\cap\mathcal{C}^{\infty}(X,L^k\otimes E).
\end{equation}
By Lichnerowicz formula \cite[(21)]{BiV}, $ 2 \overline{\partial}^{\,\ast}\overline{\partial}=\Delta^\#_k
+\frac{1}{4}K+\mathbf{c}(R)$
on  $\mathcal{C}^{\infty}(X,L^k\otimes E)$, and Corollary \ref{schrodinger} follows immediately.
This observation motivated our work.
\end{rem}
\begin{rem} As in \cite{BiV1}, we  assume that $(L, h^L, \nabla^L)$ is a positive 
Hermitian vector bundle, i.e. the curvature $R^L$ is an $\End (L)$--valued $(1,1)$--form, and for any $u\in T^{(1,0)}X 
\smallsetminus\{0\}$, $s\in L\smallsetminus\{0\}$, $\langle R^L (u,\overline{u}) s, 
\overline{s}\rangle >0$. Let $S^k (L)$ be the $k^{\text{th}}$ symmetric tensor power of $L$. Then if we replace $L^k$ in Sections \ref{s2}, \ref{s3} by $S^k(L)$, or by the irreducible representations of $L$, which are associated with the weight $ka$ (where $a$ is a given weight), when $k$ tends to $+\infty$, the techniques used in our paper still apply.
\end{rem}
\section{Covering manifolds}\label{s4}
We extend in this section our results to covering manifolds. 
\subsection{Covering manifolds, von Neumann dimension}
We present here some generalities about elliptic operators on covering manifolds and
$\g$--dimension. For details, the reader is referred to \cite[\S 4]{At}, \cite[\S{1},\S{3}]{Sh0}.

Let $\tx$ be a paracompact smooth manifold, such that there is a discrete group
$\g$ acting freely on $\tx$ having a compact quotient $X=\tx/\g$.
Let $g^{T\tx}$ be a $\g$--invariant metric on $T\tx$. Let 
$p:\tx\longrightarrow X$ be the projection.  

For a $\g$--invariant hermitian vector bundle $(\tf,h^{\tf})$, we denote by 
$\mathcal{C}_c^{\infty}(\tx,\tf)$ the space of compactly supported sections. 
Then $g^{T \tx}$, $h^{\tf}$ define an $L_2$--scalar product on $\mathcal{C}_c^{\infty}(\tx,\tf)$
as in \eqref{l2}. 
The corresponding $L_2$ space is denoted by $L_2(\tx,\tf)$.

We have a decomposition
$L_2(\tx,\tf)\cong L_2\g\otimes\mathcal{H}$ where $\mathcal{H}=L_2(U,\tf)$
is the $L_2$ space over the relatively compact
fundamental domain $U$ of the $\g$ action. This makes $L_2(\tx,\tf)$ into a 
free Hilbert $\g$--module. Since $\g$ acts by left translations $l_\gamma$ on $L_2\g$,
we obtain a unitary action of $\g$ on $L_2(\tx,\tf)$ by left translations
$L_\gamma=l_{\gamma}\otimes\Id$.
We will consider in the sequel closed $\g$--invariant subspaces of $L_2(\tx,\tf)$
for this action, called (projective) $\g$--modules.

Let  ${\mathcal A}_{\g}$ be the 
von Neumann algebra which 
consists of all bounded linear operators in $L^2\g\otimes\mathcal H$ which commute to the action of $\g$.
Let   ${\mathcal R}_{\g}$ be the von Neumann algebra of all 
bounded operators on
$L^2\g$ which commute with all $l_\gamma$. Then ${\mathcal R}_{\g}$ is generated by all right translations. Let ${\mathcal B}(\mathcal H)$  be the 
algebra  of all bounded operators on $\mathcal H$.
Then ${\mathcal A}_{\g}={\mathcal R}_{\g}\otimes {\mathcal B}(\mathcal H) $.

 If we consider 
the orthonormal basis $(\delta_\gamma)_\gamma$ in $L^2\g$, where $\delta_\gamma$ is the Dirac delta 
function at $\gamma \in \Gamma$, then the matrix of any operator $A\in{\mathcal R}_{\g}$ has the property that all
its diagonal elements are equal. Therefore we define a natural trace on ${\mathcal R}_{\g}$ as 
the diagonal element, that is,
$\operatorname{tr}_{\g} A= (A\delta_e,\delta_e)$
where $e$ is the neutral element.  
Let  $\operatorname{Tr}$ be the usual trace on ${\mathcal B}(\mathcal H)$, then we define  a trace on
${\mathcal A}_{\g}$ by $\operatorname{Tr}_{\g} = \operatorname{tr}_{\g} \otimes\operatorname{Tr}$.

For any closed $\g$--invariant space $V\subset L^2\g\otimes\mathcal H$ i.e. for any $\g$--module, the projection 
$P_V\in{\mathcal A}_{\g}$ and we define $\dim_{\g}V=\operatorname{Tr}_{\g}P_V$. In general
the $\g$--dimension is an element of $[0,\infty]$.
We also need the following fact \cite[p.398]{Sh}.
\begin{proposition}\label{inj}
Let $A:V_1\longrightarrow V_2$ be a bounded linear operator 
between two $\g$--modules, commuting with the action of $\g$. Then $\ke{A}=0$ implies
$\dim_{\g}V_1\leqslant\dim_{\g}V_2$.
\end{proposition}

Consider an elliptic, $\g$--invariant, formally self--adjoint differential 
operator $\tp$ defined in the first instance on $\mathcal{C}_c^{\infty}(\tx,\tf)$.
By a theorem of Atiyah \cite[Proposition 3.1]{At}, the minimal extension of $\tp$
(i.e. the operator closure of $\tp$) and the maximal extension of $\tp$ (i.e. $\tp^\ast$)
coincide. Hence
\begin{lemma}[Atiyah]\label{atiyah}
$\tp$ defined on $\mathcal{C}_c^{\infty}(\tx,\tf)$ is essentially self--adjoint.
\end{lemma}
Therefore $\tp$ has a unique self--adjoint extension, namely its closure. From now on, we always work
with this extension of $\tp$, which we will
denote with the same symbol. 

Then the self--adjoint extension $\tp$, as well as its
spectral projections commute with the action of $\g$. In particular, the 
spectral spaces are $\g$--modules. For a Borel set $B\subset\R$, we denote by $E(B,\tp)$ the spectral projection
corresponding to the subset $B$, and for $\mu\in \R$, set $E_\mu(\tp)=E\big((-\infty,\mu],\tp\big)$.
We introduce now a quantitative characteristic of the spectrum,
namely the von Neumann spectrum distribution function. 
For $\mu\in \R$, set
$$N_{\g}(\mu,\tp):=\tr_{\g}E_\mu(\tp)=\dim_{\g}\ran E_\mu(\tp)\,.
$$
It is non--decreasing
and the spectrum of $\tp$ coincides with the points of growth of $N_{\g}(\mu,\tp)$.
If $\tp$ is semi--bounded from below, we have $\ran E_\mu(\tp)\subset\Dom{\tp^m}$ for $m\in\N$.
Using the uniform Sobolev spaces \cite[pp. 511--2]{Sh0}, it is easily seen that $\ran E_\mu(\tp)\subset
\mathcal{C}^{\infty}(\tx,\tf)$, so that $E_\mu(\tp):L_2(\tx,\tf)\longrightarrow\mathcal{C}^{\infty}(\tx,\tf)$
is linear continuous. Let $K_\mu(\tsx,\ty)$ be the kernel of $E_\mu (\tp)$ with respect to the riemannian volume 
$dv_{\tx}$ of $g^{T \tx}$. By Schwartz kernel theorem,  
 $K_\mu(\tsx,\ty)$ is smooth. 
By \cite[Lemma 4.16]{At},  
\begin{equation*}
N_{\g}(\mu,\tp)=\tr_{\g}E_\mu(\tp)= \int_U \tr K_\mu(\tsx, \tsx)\,dv_{\tx} < + \infty.  
\end{equation*} 

\subsection{The \spin Dirac operator on a covering manifold}
Assume that there exists a $\g$--invariant pre--quantum line bundle 
$\tl$ on $\tx$  and a $\g$--invariant connection $\nabla^{\tl}$ such that 
$\tom=\frac{\sqrt{-1}}{2\pi}(\nabla^{\tl})^2$ is non--degenerate.
Let  $(\te,h^{\te})$  be a $\g$--invariant hermitian vector bundle. Let $\nabla^{\te}$ be a $\g$-invariant hermitian connection on $\te$.
Let $\tj$ be an $\g$-invariant almost complex structure on $T\tx$ such that $\tj$ is  compatible with 
$\tom$ and $g^{T \tx}$. Let $\tj_0 \in {\rm End} (TX)$ be defined by 
 \begin{equation*}
\tom (u,v) = g^{T \tx} (\tj_0 u, v), \ \ {\rm for}\ \  u,v \in T\tx.
\end{equation*}
Then $\tj$ commutes with $\tj_0$
and $\tj_0,g^{T \tx},\tom,\tj$ are the pull-back of the corresponding objects in Section \ref{s2} by $p: \tx \to X$.

We use in the sequel the same notation as in Section 2 for the corresponding objects on $X$.
Following Section \ref{s2},  we introduce the $\g$--invariant \spin Dirac operator 
$\tdk$ on $\Omega^{0,\scriptscriptstyle{\bullet}}(\tx,\tl^k\otimes\te)$
and the $\g$--invariant Laplacian 
$\tdel=\big(\nabla^{\tl^k\otimes\te}\big)^{\ast}\,\nabla^{\tl^k\otimes\te}$
on $\mathcal{C}^{\infty}(\tx,\tl^k\otimes\te)$.
Let $\tdk^{+}$ and $\tdk^{-}$  be the restrictions of $\tdk$ to 
$L_2^{0,\,\text{even}}(\tx,\tl^k\otimes\te)$ and
$L_2^{0,\,\text{odd}}(\tx,\tl^k\otimes\te)$, respectively.

\begin{proposition}
There exists $C>0$ such that for $k\in\N$, $\tdel-k\cdot\tau\circ{p}\geqslant - C$ on $L_2(\tx,\tl^k\otimes\te)$.
\end{proposition}
\begin{proof}
By applying Lichnerowicz formula \eqref{Lich} for $s\in\mathcal{C}_c^{\infty}(\tx,\tl^k\otimes\te)$,
we obtain as in the proof of Corollary \ref{bounded}, that there exists $C>0$ such that
$\big((\tdel-k\cdot\tau\circ{p})s,s\big)\geqslant-C\norm{s}^2$. 
By Lemma \ref{atiyah},
this is valid for any $s\in\Dom(\tdel-k\cdot\tau\circ{p})$.
\end{proof}
In the same vein, we can generalize Theorem 2.5.
\begin{theorem}
There exists $C>0$ such that for $k\in\N$ and any $s\in\Dom(\tdk)$
with vanishing degree zero component, 
\begin{equation}\label{main2}
\norm{\tdk{s}}^2\geqslant(2k\lambda-C)\,\norm{s}^2\,.
\end{equation}
\end{theorem}

As an immediate application of the estimate \eqref{main2} for the Dirac operator and Remark \ref{rem2.3},
we get the following asymptotic vanishing theorem which is the main result in \cite[Theorem 2.6]{Bra}.
\begin{corollary}\label{vanish1}
$\ke\tdk^{-}=\{0\}$ for large enough $k$.
\end{corollary}
We have also an analogue of Theorem \ref{specDirac}.
\begin{corollary}
There exists $C>0$ such that for $k\in\N$, the spectrum of $\tdk^2$ is contained in the
set $\{0\}\cup(2k\lambda-C,+\infty)$. 
\end{corollary}
\begin{proof}
The proof of Theorem \ref{specDirac} does not use the fact that the spectrum is discrete.
Therefore it applies in this context, too.
\end{proof}
We study now the spectrum of the $\g$--invariant Schr\"odinger operator $\tdel-k\cdot\tau\circ{p}$.
\begin{corollary}
The spectrum of the Schr\"odinger
operator $\tsc=\tdel-k\cdot\tau\circ{p}$
is contained in the union $(-a,a)\cup(2k\lambda-b,+\infty)$ ,
where $a$ and $b$ are positive constants independent of $k$.
For large enough $k$, the $\g$--dimension $d_k$ of the spectral space $E\big((-a,a),\tsc\big)$ 
corresponding to $(-a,a)$ satisfies
$d_k=\langle\ch(L^k\otimes E)\td(X),[X]\rangle$.
In particular $d_k\sim k^n(\rank{E})\vol_{\om}(X)$.
\end{corollary}
\begin{proof}
By repeating the proof of Corollary \ref{schrodinger}, we get estimates \eqref{est5} and \eqref{est11}
for smooth elements with compact support. Lemma \ref{atiyah} yields then
\begin{subequations}
\begin{align}
\big|\big(\tsc s,s\big)\big|\leqslant a\norm{s_0}^2\,,\quad s\in\Dom(\tsc)\cap P_k(\ke\tdk^{+})\,,
\label{hip2}\\
(\tsc s,s)\geqslant(2 k\lambda-b)\norm{s}^2\,,
\quad s\in\Dom(\tsc)\cap(\ke\tdk^{+})^\perp\,. \label{hip1}
\end{align}
\end{subequations}

Recall that $P_k$ represents the projection $L_2^{0,\,\scriptscriptstyle{\bullet}}(\tx,\tl^k\otimes\te)
\longrightarrow L_2^{0,0}(\tx,\tl^k\otimes\te)$. 
Since the curvatures of all our bundles are $\g$--invariant, estimate \eqref{bu} extends to the covering context with the
same proof. In particular, $P_k:\ke\tdk^{+}\longrightarrow P_k(\ke\tdk^{+})$ is bijective, 
$P_k\upharpoonright_{\ke\tdk^{+}}$ and its inverse are bounded. So 
$P_k(\ke\tdk^{+})$ is closed.
By Proposition \ref{inj},
\begin{equation}\label{dim}
\dim_{\g}\ke\tdk^{+}=\dim_{\g}P_k(\ke\tdk^{+})\,.
\end{equation}
As in \eqref{est15},  we have
\begin{equation}\label{decom}
\Dom(\tsc)=P_k(\ke\tdk^{+})\oplus(\ke\tdk^{+})^{\perp}\cap\Dom(\tsc)\,.
\end{equation}
We use now a suitable form of the minimax principle from \cite[Lemma 2.4]{Sh}:
\begin{equation}\label{var}
N_{\g}(\mu,\tsc)=\sup\{\dim_{\g}V\,:\,V\subset\Dom{\tsc}\;;\;\big(\tsc f,f\big)\leqslant\mu\norm{f}^2\;,
\forall\,f\in V\}
\end{equation}
where $V$ runs over the $\g$--modules of $L_2(\tx,\tl^k\otimes\te)$.

\noindent
By \eqref{inj}, \eqref{hip2} and  \eqref{var}, we get
\begin{equation}\label{below}
N_{\g}(a,\tsc)\geqslant\dim_{\g}\ke\tdk^{+}\,.
\end{equation}
Let us consider $\nu<2k\lambda-b$. We prove that 
\begin{equation}\label{above}
N_{\g}(\nu,\tsc)\leqslant\dim_{\g}\ke\tdk^{+}\,.
\end{equation}

Let $V\subset\Dom(\tsc)$ be an arbitrary $\g$--module with
$\big(\tsc u,u\big)\leqslant\nu\norm{u}^2$. 
If $\dim_{\g}V>\dim_{\g}\ke\tdk^{+}$, by Proposition \ref{inj} and \eqref{decom}, there exists 
$0\neq{v}\in{V}\cap(\ke D_k^{+})^\perp$, which in view of \eqref{hip1} is a contradiction.
Therefore $\dim_{\g}V\leqslant\dim_{\g}\ke\tdk^{+}$.  By \eqref{var}, we get \eqref{above}.

By  \eqref{below} and \eqref{above}, we know  that the function $N_{\g}(\nu,\tsc)$ is constant
in the interval $\nu\in[a, 2 k\lambda-b)$ and equal to $\dim_{\g}\ke\tdk^{+}$.
Enlarging a bit $a$ if necessary, we see that the spectrum of $\tsc$ is indeed contained in
$(-a,a)\cup(2 k\lambda-b,+\infty)$, and the $\g$--dimension $d_k$ of the spectral space $E\big((-a,a),\tsc\big)$
equals $\dim_{\g}\ke\tdk^{+}$. 

By Corollary \ref{vanish1},  $\dim_{\g}\ke\tdk^{+}=\ind_{\g}\tdk^{+}$. 
Moreover, Atiyah's $L_2$ index theorem
\cite[Theorem 3.8]{At} shows that $\ind_{\g}\tdk^{+}=\ind D_k^{+}$.

By \eqref{index}, the proof is achieved.
\end{proof}
\subsection*{Aknowledgements}
We thank Prof. J. M. Bismut for useful conversations. We also thank
the referee for careful reading
and  helpful comments.

\end{document}